\documentclass[12pt]{amsart}
\usepackage{latexsym,amsmath,amsfonts,amssymb,amsthm}
\theoremstyle{plain}

\newtheorem{Thm}{Theorem}
\newtheorem{Lem}{Lemma}

\theoremstyle{definition}
\newtheorem*{Ack}{Acknowledgment}
\theoremstyle{remark}
\newtheorem*{Rem}{Remark}

\def\floor #1{\left\lfloor{#1}\right\rfloor}

\begin{document}
\title{Note on some congruences of Lehmer}
\author{Hui-Qin Cao}
\address{Department of Mathematics, Nanjing University,
Nanjing 210093, People's Republic of China}
\email{caohq@nau.edu.cn}
\author{Hao Pan}
\address{Department of Mathematics, Nanjing University,
Nanjing 210093, People's Republic of China}
\email{haopan79@yahoo.com.cn}
\date{}
\maketitle

\begin{abstract}
In the paper, we generalize some congruences of Lehmer and prove that for any positive integer $n$ with $(n,6)=1$
$$
\sum_{\substack{r=1\\ (r,n)=1}}^{\floor{n/3}}\frac{1}{n-3r}\equiv\frac12
q_n(3)-\frac14 nq_n^2(3)\pmod{n^2},
$$
$$
\sum_{\substack{r=1\\ (r,n)=1}}^{\floor{n/4}}\frac{1}{n-4r}\equiv\frac34
q_n(2)-\frac38 nq_n^2(2)\pmod{n^2}
$$
and
$$
\sum_{\substack{r=1\\
(r,n)=1}}^{\floor{n/6}}\frac{1}{n-6r}\equiv\frac13 q_n(2)+\frac14
q_n(3)-n\bigg(\frac16 q_n^2(2)+\frac18 q_n^2(3)\bigg)\pmod{n^2},
$$
where $q_n(a)=(a^{\phi(n)}-1)/n$.
\end{abstract}

Let $n$ be a positive integer and let $a$ be an integer prime to $n$.
Define the Fermat quotient $q_n(a)$ by
$$
q_n(a)=\frac{a^{\phi(n)}-1}{n},
$$
where $\phi(n)$ denotes Euler's totient function. As early as
1938, Lehmer \cite{L} (or see \cite[Sections 8-9]{G}) established the following interesting
congruence:
\begin{equation}
\label{lc1}
\sum_{r=1}^{(p-1)/2}\frac{1}{r}\equiv-2q_p(2)+pq_p^2(2)\pmod{p^2}
\end{equation}
for any odd prime $p$. In 2002, Cai \cite{C} extended Lehmer's result (\ref{lc1}) and showed that
\begin{equation}
\label{cc}
\sum_{\substack{r=1\\ (r,n)=1}}^{(n-1)/2}\frac{1}{r}\equiv-2q_n(2)+nq_n^2(2)\pmod{n^2}
\end{equation}
for arbitrary odd integer $n>1$. The original proof of Cai's
congruence (\ref{cc}) depends on an identity involving generalized
Bernoulli numbers and character sums, which was proved by Szmidt,
Urbanowicz and Zagier \cite{SUZ}. However, in the same paper,
Lehmer also proved other similar congruences as follows:
\begin{equation}
\label{lc2}
\sum_{r=1}^{\floor{p/3}}\frac{1}{p-3r}\equiv\frac{1}{2}q_p(3)-\frac{1}{4}pq_p^2(3)\pmod{p^2},
\end{equation}
\begin{equation}
\label{lc3}
\sum_{r=1}^{\floor{p/4}}\frac{1}{p-4r}\equiv\frac{3}{4}q_p(2)-\frac{3}{8}pq_p^2(2)\pmod{p^2}
\end{equation}
and
\begin{equation}
\label{lc4}
\sum_{r=1}^{\floor{p/6}}\frac{1}{p-6r}\equiv\frac{1}{3}q_p(2)+\frac{1}{4}q_p(3)-p\bigg(\frac{1}{6}q_p^2(2)+\frac{1}{8}q_p^2(3)\bigg)\pmod{p^2}
\end{equation}
for each $p\geq 5$. In this note, we shall generalize Lehmer's congruences (\ref{lc2}), (\ref{lc3}) and (\ref{lc4}) for arbitrary positive integer $n$ with $(n,6)=1$.
And different from Cai's method, our approach is elementary.
\begin{Thm}
\label{mainthm}
Let $n>1$ be an integer with $(n,6)=1$. Then
\begin{equation}
\label{clc2}
\sum_{\substack{r=1\\ (r,n)=1}}^{\floor{n/3}}\frac{1}{n-3r}\equiv\frac12
q_n(3)-\frac14 nq_n^2(3)\pmod{n^2},
\end{equation}
\begin{equation}
\label{clc3} \sum_{\substack{r=1\\
(r,n)=1}}^{\floor{n/4}}\frac{1}{n-4r}\equiv\frac34 q_n(2)-\frac38
nq_n^2(2)\pmod{n^2}
\end{equation}
and
\begin{equation}
\label{clc4}
\sum_{\substack{r=1\\ (r,n)=1}}^{\floor{n/6}}\frac{1}{n-6r}\equiv\frac13
q_n(2)+\frac14 q_n(3)-n\bigg(\frac16 q_n^2(2)+\frac18
q_n^2(3)\bigg)\pmod{n^2}.
\end{equation}
\end{Thm}

\begin{Rem}
Some examples show that (\ref{clc2}) fails for $n\equiv 2, 4\pmod{6}$ and (\ref{clc3}) fails for $n\equiv 3\pmod{6}$.
\end{Rem}
\setcounter{Lem}{1}

To prove Theorem \ref{mainthm} we need some lemmas. The Bernoulli
numbers $B_n$ $(n\in\mathbb{N})$ are defined by
$$
\sum_{n=0}^\infty B_n\frac{t^n}{n!}=\frac{t}{e^t-1}.
$$
And define the Bernoulli polynomials $B_n(x)$ $(n\in\mathbb{N})$
by
$$
B_n(x)=\sum_{k=0}^n\binom{n}{k}B_{n-k}x^k.
$$
It is well-known that
\begin{equation}
\sum_{r=0}^{n-1}(x+r)^m=\frac{B_{m+1}(x+n)-B_{m+1}(x)}{m+1}.
\end{equation}
\begin{Lem}
\label{l1}
Let $p$ be a prime and $\alpha$ be a positive integer. Then
$$
\phi(p^{\alpha})\equiv
p^{\alpha}B_{\phi(p^{2\alpha})}\pmod{p^{2\alpha}}.
$$
\end{Lem}
\begin{proof}
From Euler's totient theorem, it follows that
\begin{align*}
\phi(p^{\alpha})\equiv\sum_{i=1}^{p^{\alpha}-1}i^{\phi(p^{2\alpha})}
=&\frac{1}{\phi(p^{2\alpha})+1}(B_{\phi(p^{2\alpha})+1}(p^{\alpha})-B_{\phi(p^{2\alpha})+1})\\
=&\sum_{j=1}^{\phi(p^{2\alpha})+1}\binom{\phi(p^{2\alpha})}{j-1}p^{j\alpha}B_{\phi(p^{2\alpha})+1-j}\\
\equiv&p^{\alpha}B_{\phi(p^{2\alpha})}\pmod{p^{2\alpha}},
\end{align*}
where we apply the von Staudt-Clausen theorem \cite[Theorem 5.10]{W} in the last congruence.
\end{proof}
\begin{Lem}
\label{l2} Let $p\geq 5$ be a prime and $n$ be a multiple of $p$.
Assume that $n=p^{\alpha}q$ where $p\nmid q$. Then
$$
\sum_{\substack{r=1\\ p\nmid r}}^{\floor{n/3}}\frac{1}{n-3r}
\equiv\frac{1}{2}q_{p^{2\alpha}}(3)\pmod{p^{2\alpha}},
$$
\begin{align*}
\sum_{\substack{r=1\\ p\nmid r}}^{\floor{n/4}}\frac{1}{n-4r}
\equiv\frac{3}{4}q_{p^{2\alpha}}(2)\pmod{p^{2\alpha}}
\end{align*}
and
\begin{align*}
\sum_{\substack{r=1\\ p\nmid r}}^{\floor{n/6}}\frac{1}{n-6r}
\equiv\frac{1}{3}q_{p^{2\alpha}}(2)+\frac{1}{4}q_{p^{2\alpha}}(3)\pmod{p^{2\alpha}}.
\end{align*}
\end{Lem}
\begin{proof}
Let $d\in\{3,4,6\}$ and $s$ be the least non-negative residue of
$n$ modulo $d$. Then by noting that $\phi(p^{2\alpha})-1\geq
2\alpha$, we have
\begin{align*}
\sum_{\substack{r=1\\ p\nmid
r}}^{\floor{n/d}}\frac{1}{n-dr}\equiv&\sum_{r=1}^{\floor{n/d}}(n-dr)^{\phi(p^{2\alpha})-1}
=d^{\phi(p^{2\alpha})-1}\sum_{r=1}^{\floor{n/d}}(n/d-r)^{\phi(p^{2\alpha})-1}\\
=&\frac{d^{\phi(p^{2\alpha})-1}}{\phi(p^{2\alpha})}
(B_{\phi(p^{2\alpha})}(n/d)-B_{\phi(p^{2\alpha})}(s/d))\\
=&\frac{d^{\phi(p^{2\alpha})-1}}{\phi(p^{2\alpha})}
\bigg(\sum_{i=0}^{\phi(p^{2\alpha})}\binom{\phi(p^{2\alpha})}{i}(n/d)^iB_{\phi(p^{2\alpha})-i}-B_{\phi(p^{2\alpha})}(s/d)\bigg)\\
\equiv&\frac{d^{\phi(p^{2\alpha})-1}}{\phi(p^{2\alpha})}(B_{\phi(p^{2\alpha})}-B_{\phi(p^{2\alpha})}(s/d))\pmod{p^{2\alpha}}.
\end{align*}
Since $(n,6)=1$, we have $n\equiv\pm 1\pmod{d}$, i.e., $s=1$ or
$d-1$. We know (cf. [L]) that for even integer $m>0$
$$
B_{m}(1/3)=B_{m}(2/3)=\frac{1-3^{m-1}}{2\cdot3^{m-1}}B_{m},
$$
$$
B_{m}(1/4)=B_{m}(3/4)=\frac{1-2^{m-1}}{2^{2m-1}}B_{m}
$$
and
$$
B_{m}(1/6)=B_{m}(5/6)=\frac{(1-2^{m-1})(1-3^{m-1})}{2^{m}\cdot
3^{m-1}}B_{m}.
$$
So applying Lemma \ref{l1}, we obtain that
\begin{align*}
\sum_{\substack{r=1\\ p\nmid r}}^{\floor{n/3}}\frac{1}{n-3r}
\equiv&\frac{3^{\phi(p^{2\alpha})-1}}{\phi(p^{2\alpha})}\cdot\frac{3^{\phi(p^{2\alpha})}-1}{2\cdot
3^{\phi(p^{2\alpha})-1}}B_{\phi(p^{2\alpha})}\\
\equiv&\frac{3^{\phi(p^{2\alpha})}-1}{2p^{2\alpha}}=\frac
12q_{p^{2\alpha}}(3)\pmod{p^{2\alpha}},
\end{align*}
\begin{align*}
\sum_{\substack{r=1\\ p\nmid r}}^{\floor{n/4}}\frac{1}{n-4r}
\equiv&\frac{4^{\phi(p^{2\alpha})-1}}{\phi(p^{2\alpha})}\cdot\frac{2^{2\phi(p^{2\alpha})}+2^{\phi(p^{2\alpha})}-2}{
4^{\phi(p^{2\alpha})}}B_{\phi(p^{2\alpha})}\\
\equiv&\frac{(2^{\phi(p^{2\alpha})}-1)(2^{\phi(p^{2\alpha})}+2)}{4p^{2\alpha}}
\equiv\frac{3}{4}q_{p^{2\alpha}}(2)\pmod{p^{2\alpha}},
\end{align*}
and
\begin{align*}
\sum_{\substack{r=1\\ p\nmid r}}^{\floor{n/6}}\frac{1}{n-6r}
&\equiv\frac{6^{\phi(p^{2\alpha})-1}}{\phi(p^{2\alpha})}\cdot\frac{6^{\phi(p^{2\alpha})-1}+2^{\phi(p^{2\alpha})-1}+3^{\phi(p^{2\alpha})-1}-1}{2\cdot
6^{\phi(p^{2\alpha})-1}}B_{\phi(p^{2\alpha})}\\
&\equiv\frac{1}{12p^{2\alpha}}\big(6^{\phi(p^{2\alpha})}-1+3(2^{\phi(p^{2\alpha})}-1)+2(3^{\phi(p^{2\alpha})}-1)\big)\\
&=\frac{1}{4}q_{p^{2\alpha}}(2)+\frac{1}{6}q_{p^{2\alpha}}(3)+\frac{\big((2^{\phi(p^{2\alpha})}-1)+1\big)\big((3^{\phi(p^{2\alpha})}-1)+1\big)-1}{12p^{2\alpha}}\\
&\equiv\frac{1}{3}q_{p^{2\alpha}}(2)+\frac{1}{4}q_{p^{2\alpha}}(3)\pmod{p^{2\alpha}}.
\end{align*}
\end{proof}

\begin{Lem}
\label{l3}
Let $n>1$ and $a$ be integers with $(n,6a)=1$. Then
$$
q_{n^2}(a)\equiv q_{n}(a)-\frac{1}{2}n q_{n}^2(a)\pmod{n^2}.
$$
\end{Lem}
\begin{proof}
Since $\phi(n^2)=n\phi(n)$, we have
\begin{align*}
a^{\phi(n^2)}-1=((a^{\phi(n)}-1)+1)^n-1 &\equiv
n(a^{\phi(n)}-1)+\binom{n}{2}(a^{\phi(n)}-1)^2\\
&\equiv n^2q_{n}(a)-\frac{1}{2}n^3q_{n}^2(a)\pmod{n^4},
\end{align*}
whence
$$
q_{n^2}(a)\equiv q_{n}(a)-\frac{1}{2}n q_{n}^2(a)\pmod{n^2}.
$$
\end{proof}
\begin{Lem}
\label{l4} Let $n>1$ and $a$ be integers with $(a,n)=1$. Suppose
that $n=p^{\alpha}q$ where $p$ is a prime, $\alpha>0$ and $p\nmid
q$. Then
$$
2q_{n}(a)-nq_{n}^2(a)\equiv\frac{\phi(q)}{q}\big(2q_{p^{\alpha}}(a)-p^{\alpha}q_{p^{\alpha}}^2(a)\big)\pmod{p^{2\alpha}}.
$$
\end{Lem}
\begin{proof} Compute
\begin{align*}
&2q_{n}(a)-nq_{n}^2(a)=\frac{2}{n}(a^{\phi(n)}-1)-\frac{1}{n}(a^{\phi(n)}-1)^2\\
=&\frac{2}{n}\big((a^{\phi(p^{\alpha})}-1+1)^{\phi(q)}-1\big)
-\frac{1}{n}\big((a^{\phi(p^{\alpha})}-1+1)^{\phi(q)}-1\big)^2\\
\equiv
&\frac{2}{n}\bigg(\phi(q)(a^{\phi(p^{\alpha})}-1)+\binom{\phi(q)}{2}(a^{\phi(p^{\alpha})}-1)^2\bigg)
-\frac{1}{n}\phi(q)^2(a^{\phi(p^{\alpha})}-1)^2\\
=&\frac{2\phi(q)}{n}(a^{\phi(p^{\alpha})}-1)-\frac{\phi(q)}{n}(a^{\phi(p^{\alpha})}-1)^2\\
=&\frac{\phi(q)}{q}\big(2q_{p^{\alpha}}(a)-p^{\alpha}q_{p^{\alpha}}^2(a)\big)\pmod{p^{2\alpha}}.
\end{align*}
\end{proof}
\begin{proof}[Proof of Theorem \ref{mainthm}]
 Let $p$ be an arbitrary prime factor of $n$, and assume that $n=p^{\alpha}q$ where $p\nmid q$.
Let $\mu(m)$ denote the M\"obius function. Recall (cf. \cite[Propostion 2.2.3]{IR}) that
$$
\sum_{s\mid m}\mu(s)=\begin{cases}
1\quad&\text{if }m=1,\\
0\quad&\text{if }m>1.
\end{cases}
$$
Therefore
\begin{align*}
\sum_{\substack{r=1\\ (r,n)=1}}^{\floor{n/d}}\frac{1}{n-dr}=&\sum_{\substack{r=1\\
p\nmid r,\ (r,q)=1}}^{\floor{n/d}}\frac{1}{n-dr}
=\sum_{\substack{r=1\\ p\nmid r}}^{\floor{n/d}}\frac{1}{n-dr}\sum_{s\mid(r,q)}\mu(s)\\
=&\sum_{s\mid q}\mu(s)\sum_{\substack{r=1\\ p\nmid r,\ s\mid r}}^{\floor{n/d}}\frac{1}{n-dr}
=\sum_{s\mid q}\frac{\mu(s)}{s}\sum_{\substack{t=1\\ p\nmid t}}^{\floor{n/(sd)}}\frac{1}{n/s-dt}.
\end{align*}
Observe that
$$
\sum_{s\mid q}\frac{\mu(s)}{s}=\frac{\phi(q)}{q}.
$$
Thus by Lemmas \ref{l2}, \ref{l3} and \ref{l4}, we have
\begin{align*}
&\sum_{\substack{r=1\\ (r,n)=1}}^{\floor{n/3}}\frac{1}{n-3r}
\equiv\sum_{s\mid q}\frac{\mu(s)}{s}\bigg(\frac{1}{2}q_{p^{2\alpha}}(3)\bigg)\\
\equiv&\frac{\phi(q)}{4q}\big(2q_{p^{\alpha}}(3)-p^{\alpha}q_{p^{\alpha}}^2(3)\big)
\equiv\frac{1}{4}\big(2q_n(3)-nq_n^2(3)\big)\pmod{p^{2\alpha}}.
\end{align*}
Similarly,
\begin{align*}
&\sum_{\substack{r=1\\ (r,n)=1}}^{\floor{n/4}}\frac{1}{n-4r}
\equiv\sum_{s\mid
q}\frac{\mu(s)}{s}\bigg(\frac{3}{4}q_{p^{2\alpha}}(2)\bigg)\\
\equiv&\frac{3\phi(q)}{8q}\big(2q_{p^{\alpha}}(2)-p^{\alpha}q_{p^{\alpha}}^2(2)\big)
\equiv\frac38\big(2q_n(2)-nq_n^2(2)\big)\pmod{p^{2\alpha}}.
\end{align*}
Finally,
\begin{align*}
\sum_{\substack{r=1\\ (r,n)=1}}^{\floor{n/6}}\frac{1}{n-6r}
\equiv&\sum_{s\mid q}\frac{\mu(s)}{s}\bigg(\frac{1}{3}
q_{p^{2\alpha}}(2)+\frac{1}{4}q_{p^{2\alpha}}(3)\bigg)\\
\equiv&\frac{\phi(q)}{6q}\big(2q_{p^{\alpha}}(2)-p^{\alpha}q_{p^{\alpha}}^2(2)\big)+\frac{\phi(q)}{8q}\big(2q_{p^{\alpha}}(3)-p^{\alpha}q_{p^{\alpha}}^2(3)\big)\\
\equiv&\frac{1}{6}\big(2q_n(2)-nq_n^2(2)\big)+\frac{1}{8}\big(2q_n(3)-nq_n^2(3)\big)\pmod{p^{2\alpha}}.
\end{align*}
Since $p$ is an arbitrary prime factor of $n$, we immediately deduce (\ref{clc2}), (\ref{clc3}) and (\ref{clc4})
from the Chinese remainder theorem.
\end{proof}

\begin{Ack} We thank our advisor, Professor Zhi-Wei Sun, for his helpful suggestions on this paper.
\end{Ack}


\begin{thebibliography}{9}

\bibitem {C} T.-X. Cai, \textit{A congruence involving the quotients of Euler and its applications (I)},
Acta Arith., {\bf 103}(2002), 313-320.

\bibitem {G} A. Granville, \textit{Arithmetic Properties of Binomial Coefficients I: Binomial coefficients modulo prime powers},
in Organic mathematics (Burnady,BC,1995), CMS Conf. Proc. 20, Amer. Math. Soc., Providence, RI, 1997, pp. 253-276.

\bibitem {IR} K. Ireland and M. Rosen, \textit{A classical introduction to modern number theory, 2nd ed.},
Graduate Texts in Mathematics {\bf 84}, Springer-Verlag, New York, 1990.

\bibitem {L} E. Lehmer, \textit{On congruences involving Bernoulli numbers and the quotients of Fermat and Wilson},
Ann. Math. (2), {\bf 39}(1938), 350-360.

\bibitem {SUZ} J. Szmidt, J. Urbanowicz and D. Zagier, \textit{Congruences among generalized Bernoulli numbers},
Acta Arith., {\bf 71}(1995), 273-278.

\bibitem {W} L. C. Washington, \textit{Introduction to cyclotomic fields, 2nd ed.},
Graduate Texts in Mathematics {\bf 83}, Springer-Verlag, New York, 1997.
\end{thebibliography}
\end{document}